\documentclass[twocolumn]{autart}
\usepackage{amsmath,amssymb,pstricks,color,dsfont}
\usepackage{graphicx}
\usepackage{psfrag,graphicx,epsfig}
\usepackage{algorithm,algorithmic,xspace}
\usepackage{stfloats}

\newcommand{\nnum}{\nonumber}

\newcommand{\EQ}{\begin{eqnarray}}
\newcommand{\EN}{\end{eqnarray}}
\newcommand{\EQQ}{\begin{eqnarray*}}
\newcommand{\ENN}{\end{eqnarray*}}

\newcommand{\bremark}{\begin{remark} \begin{rm} }
\newcommand{\eremark}{ \end{rm} \rule{1mm}{2mm}
\end{remark} }
\newcommand{\btheorem}{\begin{theorem} \begin{rm} }
\newcommand{\etheorem}{ \end{rm} \rule{1mm}{2mm}
\end{theorem} }
\newcommand{\blemma}{\begin{lemma} \begin{rm} }
\newcommand{\elemma}{ \end{rm} \rule{1mm}{2mm}
\end{lemma} }
\newcommand{\bcorollary}{\begin{corollary} \begin{rm} }
\newcommand{\ecorollary}{ \end{rm} \rule{1mm}{2mm}
\end{corollary} }
\newcommand{\bdefinition}{\begin{definition}\begin{rm} }
\newcommand{\edefinition}{ \end{rm} \rule{1mm}{2mm}
\end{definition} }
\newcommand{\bproposition}{\begin{proposition} \begin{rm} }
\newcommand{\eproposition}{ \end{rm} \rule{1mm}{2mm}
\end{proposition} }
\newcommand{\bexample}{\begin{example} \begin{rm} }
\newcommand{\eexample}{ \end{rm} \rule{1mm}{2mm}
\end{example} }
\newcommand{\basm}{\begin{assumption} \begin{rm}}
\newcommand{\easm}{\end{rm} 
\end{assumption}}

\newcommand{\real}{\mathds{R}}

\newcommand{\DD}{\mathcal{D}}
\newcommand{\EE}{\mathcal{E}}

\newcommand{\GG}{\mathcal{G}}
\newcommand{\HH}{\mathcal{H}}

\newcommand{\NN}{\mathcal{N}}

\newtheorem{theorem}{\bf Theorem}[section]
\newtheorem{lemma}{\bf Lemma}[section]
\newtheorem{definition}{\bf Definition}[section]
\newtheorem{remark}{\bf Remark}[section]
\newtheorem{corollary}{\bf Corollary}[section]
\newtheorem{proposition}{\bf Proposition}[section]
\newtheorem{example}{\bf Example}[section]
\newtheorem{assumption}{\bf Assumption}[section]

\newcommand\oprocendsymbol{\hbox{$\bullet$}}
\newcommand\oprocend{\relax\ifmmode\else\unskip\hfill\fi\oprocendsymbol}

\date{}

\begin{document}

\begin{frontmatter}
\title{Distributed robust adaptive equilibrium computation for generalized convex games\thanksref{footnoteinfo}}

\thanks[footnoteinfo]{This work was partially supported by Office of Naval Research grant N00014-09-1-0751. The paper was not presented at any IFAC meeting.}

\author[PSU]{Minghui Zhu\corauthref{cor}},
\corauth[cor]{Corresponding author.} \ead{muz16@psu.edu}
\author[MIT]{Emilio Frazzoli}
\ead{frazzoli@mit.edu}
\address[PSU]{Department of Electrical Engineering, Pennsylvania State University, 201 Old Main, University Park, PA, 16802, USA.}
\address[MIT]{Laboratory for Information and Decision Systems, Massachusetts Institute of Technology, 77 Massachusetts Avenue, Cambridge, MA 02139, USA.}

\maketitle

\begin{abstract} This paper considers a class of generalized convex games where each player is associated with a convex objective function, a convex inequality constraint and a convex constraint set. The players aim to compute a Nash equilibrium through communicating with neighboring players. The particular challenge we consider is that the component functions are unknown \emph{a priori} to associated players. We study two distributed computation algorithms and analyze their convergence properties in the presence of data transmission delays and dynamic changes of network topologies. The algorithm performance is verified through demand response on the IEEE 30-bus Test System. Our technical tools integrate convex analysis, variational inequalities and simultaneous perturbation stochastic approximation.
\end{abstract}
\end{frontmatter}

\section{Introduction}

Recent advances on information technologies facilitate real-time message exchanges and decision-making among geographically dispersed strategic entities. This has boosted the emergence of new generation of networked systems; e.g., the smart grid and intelligent transportation systems. These networked systems share some common features: on one hand, the entities do not belong to a single authority and may pursue different or even competitive interests; on the other hand, each entity keeps private information which is unaccessible to others. It is of great interest to design practical mechanisms which allow for efficient coordination of self-interested entities and ensure network-wide performance. Game theory along with its distributed computation algorithms represents a promising tool to achieve the goal.

In many applications, distributed computation is executed in uncertain environments. For example, mobile robots are deployed in an operating environment where environmental distribution functions are unknown to robots in advance; e.g.,~\cite{Stankovic.Johansson.Stipanovic:11,Zhu.Martinez:SICON13}. In traffic pricing, pricing policies of system operators may not be available to drivers. In optimal power flow control, the structural parameters of power systems are of national security interest and kept confidential from the public. The absence of such information makes game components; e.g., objective and constraint functions, inaccessible to players. Very recently, the informational constraint has been stimulating the investigation of \emph{adaptive} algorithms, including~\cite{Frihauf.Krstic.Basar:11,SL-MK:11,Stankovic.Johansson.Stipanovic:11} for continuous games and~\cite{JRM-HPY-GA-JSS:06,Zhu.Martinez:SICON13} for discrete games.

\emph{Literature review.} Non-cooperative game theory has been widely used as a mathematical framework to reason about multiple selfish decision makers; see for instance~\cite{Basar.Olsder:82}. These games
have found a variety of applications in economics, communication and robotics; see~\cite{Altman.Basar:98,Chung.Hollinger.Isler:11,Dockner.Jorgensen.Long.Sorger:06,Frihauf.Krstic.Basar:11,Mitchell.Bayen.Tomlin:05}. In non-cooperative games, decision making of individuals is inherently distributed. Very recently, this attractive feature has been utilized to synthesize cooperative control schemes, and a partial reference list for this regard includes~\cite{Arsie.Savla.ea:TAC09,Arslan.Marden.Shamma:07a,Li.Marden:11,Stankovic.Johansson.Stipanovic:11,Zhu.Martinez:SICON13}.

The set of papers more relevant to our work is concerned with generalized Nash games where strategy spaces are continuous and the actions of players are coupled through utility and constraint functions. Generalized Nash games are first formulated in~\cite{Arrow.Debreu:54}. Since then, a great effort has been dedicated to studying the existence and structural properties of generalized Nash equilibria in; e.g.,~\cite{Rosen:65} and the recent survey paper~\cite{Facchinei.Kanzow:07}. A number of algorithms have been proposed to compute generalized Nash equilibria, including ODE-based methods~\cite{SL-TB:87,Rosen:65}, nonlinear Gauss-Seidel-type approaches~\cite{Pang.Scutari.Facchinei.Wang:08}, iterative primal-dual Tikhonov schemes~\cite{Yin.Shanbhag.Mehta:11} and best-response dynamics~\cite{Palomar.Eldar:10}.

As mentioned, the set of papers~\cite{Frihauf.Krstic.Basar:11,SL-MK:11,JRM-HPY-GA-JSS:06,Stankovic.Johansson.Stipanovic:11,Zhu.Martinez:SICON13}
investigates the \emph{adaptiveness} of game theoretic learning algorithms. However, none of the papers mentioned in the last two paragraphs studies the \emph{robustness} of the algorithms with respect to network unreliability; e.g., data transmission delays, quantization and dynamically changing topologies. In contrast, the robustness has been extensively studied for consensus and distributed optimization, including, to name a few,~\cite{Jadbabaie.Lin.Morse:03,Nedic.Ozdaglar.Parrilo:08} for time-varying topologies,~\cite{Rabbat.Nowak:05} for quantization and~\cite{Munz.Papachristodoulou.Allgower:10} for time delays. Yet the adaptiveness issue has not been addressed in this group of papers.


\emph{Contributions.} In this paper, we aim to solve a class of generalized convex games over unreliable networks where the structures of component functions are unknown to the associated players. That is, we aim to simultaneously address the issues of adaptiveness and robustness for generalized convex games.

In the games, each player is associated with a convex objective function and subject to a private convex inequality constraint and a private convex constraint set. The component functions are assumed to be smooth and are unknown to the associated players. We investigate distributed first-order gradient-based computation algorithms for the following two scenarios:

[Scenario One] The game map is pseudo-monotone and the maximum delay (equivalently, the maximum number of packet dropouts or link breaks) is bounded but unknown;

[Scenario Two] The inequality constraints are absent, the (reduced) game map is strongly monotone and the maximum delay is known.

Inspired by simultaneous perturbation stochastic approximation for optimization in~\cite{JS:03}, we utilize finite differences with diminishing approximation errors to estimate first-order gradients. We propose two distributed algorithms for the two scenarios and formally prove their asymptotic convergence. The comparison of the two proposed algorithms is given in Section~\ref{sec:comparison}. The analysis integrates the tools from convex analysis, variational inequalities and simultaneous perturbation stochastic approximation. The algorithm performance is verified through demand response on the IEEE 30-bus Test System. A preliminary version of the current paper was published in~\cite{Zhu.Frazzoli:CDC12} where the adaptiveness issue was not investigated.


\section{Problem formulation}\label{sec:formulation}

In this section, we present the generalized convex game considered in the paper. It is followed by the notions and notations used throughout the paper.

\subsection{Generalized convex game}

Consider the set of players $V \triangleq \{1,\cdots,N\}$ where the state of player~$i$ is denoted as $x^{[i]}\in X_i\subseteq \real^{n_i}$. The players are selfish and pursue different interests. In particular, given the joint state $x^{[-i]}\in X_{-i} \triangleq \prod_{j\neq i}X_j$ of its rivals\footnote{We use the shorthand $-i\triangleq V\setminus\{i\}$ throughout the paper.}, each player~$i$ aims to solve the following program parameterized by $x^{[-i]}\in X_{-i}$: \begin{align}\min_{x^{[i]}\in X_i}f_i(x^{[i]},x^{[-i]}),\quad {\rm s.t.} \quad G^{[i]}(x^{[i]},x^{[-i]})\leq0,\label{e10}\end{align} where $f_i : \real^n\rightarrow \real$ and $G^{[i]} : \real^n\rightarrow\real^{m_i}$ with $n\triangleq \sum_{i\in V}n_i$. In the remainder of the paper, we assume that the following properties about problem~\eqref{e10} hold:

\begin{assumption} The maps  $f_i$ and $G^{[i]}$ are smooth, and the maps  $f_i(\cdot,x^{[-i]})$ and $G^{[i]}(\cdot,x^{[-i]})$ are convex in $x^{[i]}$. The set $X_i$ is convex and compact, and $X\cap Y\neq\emptyset$ where $X \triangleq \prod_{i\in V}X_i$ and $Y \triangleq \prod_{i\in V}Y_i$ with $Y_i \triangleq \{x\in X\; |\; G^{[i]}(x) \leq 0\}$.\label{asm7}
\end{assumption}


We now proceed to provide an equivalent form of problem~\eqref{e10}. To achieve this, we define the set-valued map $X_i^f : X_{-i} \rightarrow 2^{X_i}$ as follows: \begin{align*}X_i^f(x^{[-i]}) = \{x^{[i]}\in X_i \; | \; G^{[i]}(x^{[i]},x^{[-i]})\leq 0\}.\end{align*} The set $X_i^f(x^{[-i]})$ represents the collection of feasible actions for player~$i$ when its opponents choose the joint state of $x^{[-i]}\in X_{-i}$. With the map $X_i^f$, problem~\eqref{e10} of player~$i$ is equivalent to the following one: \begin{align}\min_{x^{[i]}\in X_i^f(x^{[-i]})}f_i(x^{[i]},x^{[-i]}).\label{e11}\end{align}

Given $x^{[-i]}\in X_{-i}$, each player~$i$ aims to solve problem~\eqref{e11}. The collection of such coupled optimization problems consists of the \emph{generalized convex game} (for short, CVX). For the CVX game, we adopt the \emph{generalized Nash equilibrium} (for short, GNE) as the solution notion which none of the players is willing to unilaterally deviate from:

\begin{definition} The joint state $\tilde{x}\in X\cap Y$ is a generalized Nash equilibrium of the CVX game if the following holds: \begin{align}&f_i(\tilde{x})\leq f_i(x^{[i]},\tilde{x}^{[-i]}),\quad \forall x^{[i]}\in X_i^f(\tilde{x}^{[-i]}),\quad \forall i\in V.\nnum\end{align}\label{def1}
\end{definition}

Denote by $\mathbb{X}_{\rm CVX}$ the set of GNEs of the CVX game. The following lemma verifies the non-emptiness of $\mathbb{X}_{\rm CVX}$.

\begin{lemma} The set of generalized Nash equilibria of the CVX game is not empty, i.e., $\mathbb{X}_{\rm CVX} \neq \emptyset$.\label{lem3}
\end{lemma}

\textbf{Proof:} Recall that $f_i$ is convex and $X\cap Y$ is compact. Hence, $\mathbb{X}_{\rm CVX} \neq \emptyset$ is a direct result of~\cite{Facchinei.Kanzow:07,Rosen:65}. \oprocend


In the CVX game, the players desire to seek a GNE. It is noted that $f_i$, $G^{[i]}$, and $X_i$ are private information of player~$i$ and unaccessible to others. In order to compute a GNE, it becomes necessary that the players are inter-connected and able to communicate with each other to exchange their partial estimates of GNEs. The interconnection between players will be represented by a directed graph ${\GG} = (V,\EE)$ where
$\EE\subset V\times V\setminus {\rm diag}(V)$ is the set of
edges. The neighbor relation is determined by the dependency of $f_i$ and/or $G^{[i]}$ on $x^{[j]}$. In particular, $(i,j)\in\EE$ if and only if $f_i$ and/or $G^{[i]}$ depend upon $x^{[j]}$. Denote by $\NN_i^{\rm IN} \triangleq \{j\in V\;|\; (i,j)\in\EE\}$ the set of in-neighbors of
player~$i$.

In this paper, we aim to develop distributed algorithms which allow for the computation of GNEs in the presence of the following two challenges.
\begin{enumerate}
\item Data transmissions between the players in $V$ are unreliable, and subject to transmission delays, packet dropouts and/or link breaks. We let $x^{[j]}(k - \tau^{[i]}_j(k))$ with $\tau^{[i]}_j(k)\geq0$ be either $(i)$ the outdated state of player~$j$ received by player~$i$ at time~$k$ due to transmission delays; or $(ii)$ the latest state of player~$j$ received by player~$i$ by time~$k$ due to packet dropouts and/or link breaks. For example, if the link from player $j$ to player $i$ is broken at time $k$, then player~$i$ cannot receive the message $x^{[j]}(k)$ from player $j$. For this case, player~$i$ uses the more recent received message sent from player $j$ which is $x^{[j]}(k-\tau_j^{[i]}(k))$. Denote by $\Lambda_i(k)\triangleq \{x^{[j]}(k-\tau^{[i]}_j(k))\}_{j\in\NN^{\rm IN}_i}$ the set of latest states of $\NN_i^{\rm IN}$ received by player~$i$ at time~$k$. Let $\tau_{\max} \triangleq \sup_{k\geq0}\max_{i\in V}\max_{j\in\NN_i^{\rm IN}}\tau^{[i]}_j(k)$ be the maximum delay or the maximum number of successive packet dropouts.
\item Each player~$i$ is unaware of the structures of $f_i$ and $G^{[i]}$ but can observe their realizations. That is, if the players input the joint state $x$ into $f_i$, then player~$i$ can observe the realized value of $f_i(x)$.
\end{enumerate}

The disclosure of the value $f_i(x)$ is case dependent. In~\cite{Stankovic.Johansson.Stipanovic:11,Zhu.Martinez:SICON13}, mobile sensors are unaware of environmental distribution functions but they can observe induced utilities via on-site measurements. This is an example of engineering systems. Section~\ref{sec:simulations} will provide a concrete example of demand response of power networks where the system operator discloses realized values via communication given the decisions of end-users. This is an example of social systems.

\section{Preliminaries}\label{sec:preliminaries}

In the CVX game, each player is subject to an inequality constraint. In optimization literature, Lagrangian relaxation is a widely used approach to handle inequality constraints. Following this vein, we will perform Lagrangian relaxation on the CVX game and obtain the unconstrained convex (for short, UC) game.

In the UC game, there are two sets of players: the set of primal players $V$ and the set of dual players $V_m \triangleq \{1,\cdots,N\}$. Define the following private Lagrangian $\mathcal{H}_i : \real^n\times\real^{m_i}_{\geq0}\rightarrow\real$: $\mathcal{H}_i(x,\mu^{[i]}) \triangleq f_i(x) + \langle \mu^{[i]}, G^{[i]}(x)\rangle$, for primal player~$i\in V$ and dual player~$i\in V_m$. Given any $x^{[-i]}\in X_{-i}$ and $\mu^{[i]}\in\real^{m_i}_{\geq0}$, each primal player~$i\in V$ aims to minimize $\mathcal{H}_i$ over $x^{[i]}\in X_i$; i.e., $\min_{x^{[i]}\in X_i}\mathcal{H}_i(x^{[i]},x^{[-i]},\mu^{[i]})$, and, instead, the objective of the dual player $i\in V_m$ is to maximize $\mathcal{H}_i$ over $\mu^{[i]}\in M_i\subseteq\real^{m_i}_{\geq0}$; i.e., $\min_{\mu^{[i]}\in M_i}-\mathcal{H}_i(x,\mu^{[i]})$ where the set $M_i\subseteq\real^{m_i}_{\geq0}$ is convex, non-empty and will be introduced in the sequel. We let $\eta \triangleq (x,\mu)$ and the set $K \triangleq X\times M$ with $M \triangleq \prod_{i\in V}M_i$. The above game among the players in $V\cup V_m$ is referred to as the UC game parameterized by the set $M$. The set of $M$ will play an important role in determining the properties of the UC game and we will discuss the choice of $M$ later. The solution concept for the UC game parameterized by $M$ is the standard notion of Nash equilibrium given below:

\begin{definition} The joint state of $(\tilde{x},\tilde{\mu}) \in X\times M$ is a Nash equilibrium of the UC game parameterized by $M$ if the following holds for all $i\in V$: $\mathcal{H}_i(\tilde{x},\tilde{\mu}^{[i]})\leq \mathcal{H}_i(x^{[i]},\tilde{x}^{[-i]},\tilde{\mu}^{[i]}),\quad \forall x^{[i]}\in X_i$, and the following holds for all $i\in V_m$: $\mathcal{H}_i(\tilde{x},\mu^{[i]})\leq \mathcal{H}_i(\tilde{x},\tilde{\mu}^{[i]}),\quad \forall \mu^{[i]}\in M_i$.\label{def2}
\end{definition} Denote by $\mathbb{X}_{\rm UC}(M)$ the set of NEs of the UC game parameterized by $M$. We now proceed to illustrate the relation between the UC game and the CVX game. Before doing so, let us state the following boundedness assumption on the dual solutions:

\begin{assumption} There is a vector $\vartheta \triangleq (\vartheta_i)_{i\in V}\in\real^N_{>0}$ such that for any $(\tilde{x},\tilde{\mu})\in \mathbb{X}_{\rm UC}(\real^m_{\geq0})$,  $\|\tilde{\mu}^{[i]}\| \leq \vartheta_i$ for $i\in V$.\label{asm1}
\end{assumption}

\begin{remark} We would like to make a remark on Assumption~\ref{asm1}. For convex optimization, it is shown in~\cite{Hiriart.Lemarechal:96} that the Lagrangian multipliers are uniformly bounded under the standard Slater's condition. This boundedness property is used in~\cite{Nedic.Ozdagalr:08} and further in~\cite{Zhu.Martinez:09tac,Zhu.Martinez:10tac} to solve convex and non-convex programs in centralized and distributed manners. In this paper, Assumption~\ref{asm1}, on one hand, ensures the existence of GNEs, and on the other hand, guarantees the boundedness of estimates and gradients for the case of unknown $\tau_{\max}$. We will discuss the verification of Assumption~\ref{asm1} in Section~\ref{sec:discussion}.\oprocend\label{rem6}
\end{remark}

The following proposition characterizes the relations between the UC game and the CVX game.

\begin{proposition} The following properties hold:
\begin{enumerate}
\item[(P1)] Consider any $(\tilde{x}, \tilde{\mu})\in\mathbb{X}_{\rm UC}(M)$. We have $\tilde{x}\in \mathbb{X}_{\rm CVX}$ if the following properties hold:
\begin{itemize}
\item  feasibility; i.e., $G^{[i]}(\tilde{x})\leq0$ for all $i\in V$;
\item  slackness complementarity; i.e., $\langle \tilde{\mu}^{[i]}, G^{[i]}(\tilde{x})\rangle = 0$ for all $i\in V$.
\end{itemize}
\item[(P2)] Suppose Assumption~\ref{asm1} holds. Consider the UC game parameterized by $M$ with $M_i \triangleq \{\mu^{[i]}\in\real^{m_i}_{\geq0}\;|\;\|\mu^{[i]}\|\leq \vartheta_i+r_i\}$ for some $r_i>0$. Then $\mathbb{X}_{\rm UC}(M)\neq\emptyset$. In addition, for any $(\tilde{x}, \tilde{\mu})\in\mathbb{X}_{\rm UC}(M)$, it holds that $\tilde{x}\in \mathbb{X}_{\rm CVX}$.
\end{enumerate} \label{pro5}
\end{proposition}

\textbf{Proof:} The proof of (P1) is a slight extension of the results in~\cite{Bertsekas:09} to our game setup. For the sake of completeness, we summarize the analysis here. Since $\langle \tilde{\mu}^{[i]}, G^{[i]}(\tilde{x})\rangle = 0$, we have the following relation: \begin{align}f_i(\tilde{x}) = \HH_i(\tilde{x},\tilde{\mu}^{[i]}) - \langle\tilde{\mu}^{[i]},G^{[i]}(\tilde{x})\rangle = \HH_i(\tilde{x},\tilde{\mu}^{[i]}).\label{e7}\end{align} Since $(\tilde{x},\tilde{\mu})\in\mathbb{X}_{\rm UC}(M)$, then the following relation holds for all $x^{[i]}\in X_i$: \begin{align}\mathcal{H}_i(\tilde{x},\tilde{\mu}^{[i]})
-\mathcal{H}_i(x^{[i]},\tilde{x}^{[-i]},\tilde{\mu}^{[i]}) \leq 0.\label{e17}\end{align} Substitute~\eqref{e7} into~\eqref{e17}, and it renders the following for any $x^{[i]}\in X_i$:
\begin{align}f_i(\tilde{x}) \leq f_i(x^{[i]},\tilde{x}^{[-i]}) + \langle\tilde{\mu}^{[i]},G^{[i]}(x^{[i]},\tilde{x}^{[-i]})\rangle,\nnum\end{align} and thus $f_i(\tilde{x}) \leq f_i(x^{[i]},\tilde{x}^{[-i]}),\quad \forall x^{[i]}\in X_i^f(\tilde{x}^{[-i]})$, by noting that $\langle\tilde{\mu}^{[i]},G^{[i]}(x^{[i]},\tilde{x}^{[-i]})\rangle \leq 0$ for any $ x^{[i]}\in X_i^f(\tilde{x}^{[-i]})$. Recall that $G^{[i]}(\tilde{x})\leq0$ for all $i\in V$. The above arguments hold for all $i\in V$, and thus it establishes that $\tilde{x}\in\mathbb{X}_{\rm CVX}$.

We now proceed to show (P2). Since $\mu^{[i]}\geq0$, $\HH_i$ is a positive combination of convex functions of $x^{[i]}$. Thus $\HH_i$ is convex in $x^{[i]}$. It is easy to see that $-\HH_i$ is convex (actually affine) in $x^{[i]}$. Since $X$ and $M$ are convex and compact, it follows that $\mathbb{X}_{\rm UC}(M) \neq \emptyset$ by the results on generalized convex games in; e.g.,~\cite{Facchinei.Kanzow:07,Rosen:65}. Pick any $(\tilde{x},\tilde{\mu})\in\mathbb{X}_{\rm UC}(M)$ and then $\|\tilde{\mu}^{[i]}\| \leq \vartheta_i$ by Assumption~\ref{asm1}. We then have the following relation: \begin{align}{\HH}_i(\tilde{x},\tilde{\mu}^{[i]})
  -{\HH}_i(\tilde{x},\mu^{[i]}) \geq 0,\quad \forall\mu^{[i]}\in M_i,\nnum\end{align} and thus,
\begin{align}\langle \mu^{[i]}-\tilde{\mu}^{[i]}, G^{[i]}(\tilde{x})\rangle \leq 0,\quad \forall\mu^{[i]}\in M_i.\label{e6}\end{align}

We now proceed to verify by contradiction the feasibility of $G^{[i]}(\tilde{x}) \leq 0$ and the complementary slackness of $\langle  \tilde{\mu}^{[i]}, G^{[i]}(\tilde{x})\rangle = 0$. Assume that $G^{[i]}(\tilde{x})_{\ell} > 0$. We choose $\mu^{[i]}$ such that $\mu^{[i]}_{\ell} = \tilde{\mu}^{[i]}_{\ell} + \pi r_i$ and $\mu^{[i]}_{\ell'} = 0$ for $\ell' \neq \ell$ where $\pi>0$ is sufficiently small such that $\mu^{[i]}\in M_i$. Then it follows from~\eqref{e6} that $r_i G^{[i]}(\tilde{x})_{\ell} \leq 0$ which is a contradiction. Hence, we have the feasibility of $G^{[i]}(\tilde{x})\leq 0$. Combine it with $\tilde{\mu}^{[i]}\geq0$, and it renders that $\langle \tilde{\mu}^{[i]}\rangle, G^{[i]}(\tilde{x})\leq 0$. On the other hand, we let $\mu^{[i]} = 0$ in the relation~\eqref{e6} and have $\langle \tilde{\mu}^{[i]}, G^{[i]}(\tilde{x})\rangle \geq 0$. Hence, it renders that $\langle \tilde{\mu}^{[i]}, G^{[i]}(\tilde{x})\rangle = 0$. We reach the desired result by using the property (P1). \oprocend

\subsection{Notations and notions}\label{subsection:notations}

The vector $\textbf{1}_n$ represents the column vector with $n$ ones. The vector $e^{[\ell]}_{n}$ is the one in $\real^n$ whose $\ell$-th coordinate is one and whose other coordinates are all zero. For any pair of vectors $a,b\in\real^p$, the relation $a\leq b$ means $a_{\ell}\leq b_{\ell}$ for all $1\leq\ell\leq p$. Since the function $f_i$ is continuous and $X_i$ is compact, then the following quantities are well-defined:
$\sigma_{i,\min} = \inf_{x\in X_i}f_i(x),\quad \sigma_{i,\max} = \sup_{x\in X_i}f_i(x),\quad
\sigma_{\min} = \inf_{x\in X}\sum_{i\in V}f_i(x),\quad \sigma_{\max} = \sup_{x\in X}\sum_{i\in V}f_i(x)$. Given a non-negative scalar sequence $\{\alpha(k)\}_{k\geq0}$, it is \emph{summable} if $\sum_{k=0}^{\infty}\alpha(k)<+\infty$ and \emph{square summable} if $\sum_{k=0}^{\infty}\alpha(k)^2<+\infty$.

In the remainder of the paper, we will use some notions about \emph{monotonicity}, and the readers are referred to~\cite{Aubin.Frankowska:09,Facchinei.Pang:03} for detailed discussion. The mapping $F : Z \rightarrow Z'$ is \emph{strongly monotone} with constant $\rho > 0$ over $Z$ if for each pair of $\eta,\eta'\in Z$, the following holds: \begin{align*}\langle F(\eta) - F(\eta'), \eta - \eta' \rangle \geq \rho \|\eta - \eta'\|^2.\end{align*} The mapping $F : Z \rightarrow Z'$ is \emph{monotone} over $Z$ if for each pair of $\eta,\eta'\in Z$, the following holds: \begin{align*}\langle F(\eta) - F(\eta'), \eta - \eta' \rangle \geq 0.\end{align*} The mapping $F : Z \rightarrow Z'$ is \emph{pseudo-monotone} over $Z$ if for each pair of $\eta,\eta'\in Z$, it holds that $\langle F(\eta'), \eta - \eta' \rangle \geq 0$ implies $\langle F(\eta), \eta - \eta' \rangle \geq 0$. It is known that strong monotonicity implies monotonicity, and monotonicity implies pseudo-monotonicity~\cite{Facchinei.Pang:03}.

Given the non-empty, convex, and closed set $Z\in {\real}^n$, the
projection operator onto $Z$, $P_Z: \real^n \rightarrow Z$, is
defined as $P_Z[z] = {\rm argmin}_{x\in Z}\|x-z\|$. The following is the non-expansiveness of the projection operator; e.g.,~\cite{Bertsekas:09}.

\begin{lemma} Let $Z$ be a non-empty, closed and convex set in ${\real}^n$. For any $z\in{\real}^n$, the following holds for any $y\in Z$: $\|P_X[z] - y
  \|^2 \leq \| z - y \|^2 - \| P_X[z] - z\|^2$. \label{lem6}
\end{lemma}

We define $\nabla_{x^{[i]}}\HH_i(x,\mu^{[i]})$ as the gradient of the convex function $\HH_i(\cdot,x^{[-i]},\mu^{[i]})$ at $x^{[i]}$, and $\nabla_{\mu^{[i]}}\HH_i(x,\mu^{[i]}) = G^{[i]}(x)$ as the gradient of the concave (actually affine) function $\HH_i(x,\cdot)$ at $\mu^{[i]}$. Define $\nabla \Theta$ as the map of partial gradients of the players' objective functions: \begin{align*}\nabla \Theta(\eta) &\triangleq [\nabla_{x^{[1]}}\HH_1(x,\mu^{[1]})^T\cdots \nabla_{x^{[N]}}\HH_N(x,\mu^{[N]})^T\\ &\nabla_{\mu^{[1]}}-\HH_1(x,\mu^{[1]})\cdots \nabla_{\mu^{[N]}}-\HH_N(x,\mu^{[N]})]^T.\end{align*} The map $\nabla \Theta$ is referred to as the \emph{game map}. It is well-known that $\nabla\Theta$ is monotone (resp. strongly monotone) over $K$ if and only if its Jacobian $\nabla^2\Theta(\eta)$ is positive semi-definite (resp. positive definite) for all $\eta\in K$.

For $G^{[i]}(x)\leq0$ being absent, we define $\nabla_{x^{[i]}}f_i(x)$ as the gradient of the convex function $f_i(\cdot,x^{[-i]})$ at $x^{[i]}$. Define $\nabla \Theta^r$ as the map of partial gradients of the players' objective functions: \begin{align*}\nabla \Theta^r(x) \triangleq [\nabla_{x^{[1]}}f_1(x)^T\cdots \nabla_{x^{[N]}}f_N(x)^T]^T.\end{align*} The map $\nabla \Theta^r$ is referred to as the \emph{(reduced) game map} for $G^{[i]}(x)\leq0$ being absent.


Given any convex and closed set $Z$, let $\mathbb{P}_Z$ be the projection operator onto the set $Z$.

The following lemma is a direct result of the smoothness of the component functions and the compactness of the constraint sets.

\begin{lemma} The image set of the game map $\nabla\Theta$ is uniformly bounded over $K$. In addition, the game map $\nabla\Theta$ (resp. $\nabla\Theta^r$) is Lipschitz continuous over $K$ with constant $L_{\Theta}$ (resp. $L_{\Theta^r}$). \label{lem4}
\end{lemma}

\section{Distributed computation algorithms}\label{sec:algorithm}

In this section, we will study the scenarios mentioned in the introduction. Proposition~\ref{pro5} reveals that it suffices to compute a Nash equilibrium in $\mathbb{X}_{\rm UC}$ if Assumption~\ref{asm1} holds. In the remainder of this section, we will suppose Assumption~\ref{asm1} holds, and then each dual player~$i$ can define the set of $M_i \triangleq \{\mu^{[i]}\in\real^{m_i}_{\geq0}\;|\;\|\mu^{[i]}\|\leq \vartheta_i+r_i\}$ with $r_i>0$. Assumption~\ref{asm1} will be justified in Section~\ref{sec:assumption}.

\subsection{Scenario one: pseudo-monotone game map and unknown $\tau_{\max}$}

In this section, we synthesize a distributed first-order algorithm for the case where the quantity $\tau_{\max}$ is unknown and the game map $\nabla \Theta$ is merely pseudo-monotone.

Algorithm~\ref{ta:algo1} is based on projected primal-dual gradient methods where each primal or dual player~$i$\footnote{In practical implementation, the update rules of primal player~$i$ and dual player~$i$ are both executed by real player~$i$.} updates its own estimate by moving along its partial gradient with a certain step-size and projecting the estimate onto the local constraint set. Recall that $f_i$ and $G^{[i]}$ are unknown to player~$i$. Then player~$i$ cannot compute partial gradient $\nabla_{x^{[i]}}\HH_i$. Inspired by simultaneous perturbation stochastic approximation in; e.g.,~\cite{JS:03}, we use finite differences to approximate the partial gradient $\nabla_{x^{[i]}}\HH_i$; i.e., for $\ell = 1,\cdots,n_i$, \begin{align}&[\nabla_{x^{[i]}}\HH_i(x(k),\mu^{[i]}(k))]_{\ell}\nnum\\
&\approx \frac{1}{2c_i(k)}\big(\HH_i(x^{[i]}(k)+c_i(k)e^{[\ell]}_{n_i},x^{[-i]}(k),\mu^{[i]}(k))\nnum\\
&-\HH_i(x^{[i]}(k)-c_i(k)e^{[\ell]}_{n_i},x^{[-i]}(k),\mu^{[i]}(k))\big),\label{e51}\end{align}
where $-c_i(k)e^{[\ell]}_{n_i}$ and $c_i(k)e^{[\ell]}_{n_i}$ are two-way perturbations. In addition, $\nabla_{\mu^{[i]}}\HH_i(x(k),\mu^{[i]}(k)) = G^{[i]}(x(k))$. In Algorithm~\ref{ta:algo1}, $\DD^{[i]}_x(k)$ (resp. $\DD^{[i]}_{\mu}(k)$) is the right-hand side of~\eqref{e51} (resp. $\nabla_{\mu^{[i]}}\HH_i(x,\mu^{[i]})$) with delayed estimates. Here we assume that player~$i$ can observe the \emph{values} associated with $f_i$ and $G^{[i]}$ and thus $\DD^{[i]}_x(k)$ and $\DD^{[i]}_{\mu}(k)$.

The scalar $c_i(k)>0$ is the perturbation magnitude. When it is small, the finite-difference approximation is close to the partial gradient. In order to asymptotically eliminate the error induced by the finite-difference approximation, $c_i(k)$ needs to be diminishing. The decreasing rate of $c_i(k)$ should match that of computation step-sizes $\alpha(k)$. Otherwise, the convergence to Nash equilibrium may be prevented. This is captured by (A4) of Theorem~\ref{the1}.

\begin{algorithm}[htbp]\caption{Distributed gradient-based algorithm for Scenario one} \label{ta:algo1}
\textbf{Require:} Each primal player $i\in V$ chooses the initial state $x^{[i]}(0)\in X_i$. And each dual player~$i\in V_m$ defines the set $M_i$ and chooses the initial state $\mu^{[i]}(0)\in M_i$.

\textbf{Ensure:} At each $k \ge 0$, each player in $V\cup V_m$ executes the following steps:

1. Each primal player $i\in V$ updates its state according to the following rule:\begin{align*} x^{[i]}(k+1) = \mathbb{P}_{X_i}[x^{[i]}(k) - \alpha(k)\DD^{[i]}_x(k)],
\end{align*} where the approximate gradient $\DD^{[i]}_x(k)$ is given by \begin{align}&[\DD^{[i]}_x(k)]_{\ell} = \frac{1}{2c_i(k)}
\{f_i(x^{[i]}(k)+c_i(k)e^{[\ell]}_{n_i},\Lambda_i(k))\nnum\\
&+ \langle\mu^{[i]}(k),G^{[i]}(x^{[i]}(k)+c_i(k)e^{[\ell]}_{n_i},\Lambda_i(k))\rangle\nnum\\
&-f_i(x^{[i]}(k)-c_i(k)e^{[\ell]}_{n_i},\Lambda_i(k))\nnum\\
&- \langle\mu^{[i]}(k),G^{[i]}(x^{[i]}(k)
-c_i(k)e^{[\ell]}_{n_i},\Lambda_i(k))\rangle\},\label{e38}\end{align} for $\ell = 1,\cdots,n_i$.

2. Each dual player $i\in V_m$ updates its state according to the following rule:\begin{align*}
\mu^{[i]}(k+1) = \mathbb{P}_{M_i}[\mu^{[i]}(k) + \alpha(k)\DD^{[i]}_{\mu}(k)],
\end{align*} where the gradient $\DD^{[i]}_{\mu}(k)$ is given by $\DD^{[i]}_{\mu}(k) = \nabla_{\mu^{[i]}}\HH_i(x^{[i]}(k),\Lambda_i(k),\mu^{[i]}(k)) = G^{[i]}(x^{[i]}(k),\Lambda_i(k))$.

3. Repeat for $k = k+1$.
\end{algorithm}

The following theorem demonstrates that Algorithm~\ref{ta:algo1} is able to achieve a GNE from any initial state in $K$.

\begin{theorem} Suppose the following hold:
\begin{enumerate}
\item[(A1)] the quantity $\tau_{\max}$ is finite;
\item[(A2)] Assumptions~\ref{asm7} and~\ref{asm1} hold;
\item[(A3)] the sequence of $\{\alpha(k)\}$ is positive, not summable but square summable;
\item[(A4)] the sequence of $\displaystyle{\{\alpha(k)\max_{i\in V}c_i(k)\}}$ is summable;
\item[(A5)] the game map $\nabla\Theta$ is pseudo-monotone over $K$.
\end{enumerate}
For any initial state $\eta(0)\in K$, the sequence of $\{\eta(k)\}$ generated by Algorithm~\ref{ta:algo1} converges to some $(\tilde{x},\tilde{\mu})\in{\mathbb{X}}_{\rm UC}$ where $\tilde{x}\in\mathbb{X}_{\rm CVX}$.\label{the1}
\end{theorem}

\textbf{Proof:} First of all, it is noted that $K$ is compact since each $X_i$ and $X_i$ is compact in Assumption~\ref{asm1} and (P2) in Proposition~\ref{pro5} as well as $\nabla\Theta$ is uniformly bounded over $K$ in Lemma~\ref{lem4}.

Then we write Algorithm~\ref{ta:algo1} in the following compact form:
\begin{align}\eta(k+1) = \mathbb{P}_K[\eta(k)-\alpha(k)\DD(k)],\label{e19}
\end{align} where $\DD(k) \triangleq ((\DD^{[i]}_x(k)^T)^T_{i\in V},(-\DD^{[i]}_{\mu}(k)^T)^T_{i\in V_m})^T$ is subject to time delays and perturbations. Then pick any $\eta\in K$. It follows from the non-expansiveness of the projection operator $\mathbb{P}_K$, Lemma~\ref{lem6}, that for any $\eta\in K$, the following relation holds:
\begin{align}&\|\eta(k+1)-\eta\|^2\nnum\\
&\leq \|\eta(k)-\alpha(k)\DD(k)-\eta\|^2-\|\DD(k)\|^2\nnum\\
&\leq \|\eta(k)-\alpha(k)\DD(k)-\eta\|^2\nnum\\
&= \|\eta(k)-\eta\|^2 - 2\alpha(k)\langle \DD(k), \eta(k)-\eta\rangle\nnum\\
&+ \alpha(k)^2\|\DD(k)\|^2.\label{e25}
\end{align}

It follows from~\eqref{e25} that
\begin{align}&2\alpha(k)\langle \DD(k), \eta(k)-\eta\rangle\leq \|\eta(k) - \eta\|^2\nnum\\
&- \|\eta(k+1) - \eta\|^2 +\alpha(k)^2\|\DD(k)\|^2.\label{e23}\end{align}

For any $k\geq0$ and $i$, we define the following: \begin{align}&\hat{\DD}^{[i]}_x(k) = \nabla_{x^{[i]}}\HH_i(x^{[i]}(k),\{x^{[j]}(k)\}_{j\in\NN^{\rm IN}_i},\mu^{[i]}(k)),\nnum\\
&\hat{\DD}^{[i]}_{\mu}(k) \triangleq \nabla_{\mu^{[i]}}\HH_i(x^{[i]}(k),\{x^{[j]}(k)\}_{j\in\NN^{\rm IN}_i},\mu^{[i]}(k)),\nnum\\
&\hat{\DD}(k) \triangleq ((\hat{\DD}^{[i]}_x(k)^T)^T_{i\in V},(-\hat{\DD}^{[i]}_{\mu}(k)^T)^T_{i\in V_m})^T,\label{e16}\end{align} which are the gradients evaluated at the delay-free states. So, the quantity $\hat{\DD}(k)$ is free of time delays and perturbations.

Similarly, for any $k\geq0$ and $i$, we define the following: \begin{align}&\tilde{\DD}^{[i]}_x(k) = \nabla_{x^{[i]}}\HH_i(x^{[i]}(k),\Lambda_i(k),\mu^{[i]}(k)),\nnum\\
&\tilde{\DD}^{[i]}_{\mu}(k) \triangleq \nabla_{\mu^{[i]}}\HH_i(x^{[i]}(k),\Lambda_i(k),\mu^{[i]}(k)),\nnum\\
&\tilde{\DD}(k) \triangleq ((\tilde{\DD}^{[i]}_x(k)^T)^T_{i\in V},(-\tilde{\DD}^{[i]}_{\mu}(k)^T)^T_{i\in V_m})^T,\nnum\end{align} which are the gradients evaluated at the delayed states. Then the quantity $\tilde{\DD}(k)$ is free of perturbations but subject to time delays.

With the above notations at hand, the relation~\eqref{e23} implies the following: \begin{align}&2\alpha(k)\langle \hat{\DD}(k), \eta(k)-\eta\rangle \leq \alpha(k)^2\|\DD(k)\|^2 + \|\eta(k) - \eta\|^2\nnum\\&- \|\eta(k+1) - \eta\|^2 + 2\alpha(k)\langle\tilde{\DD}(k)-\DD(k),\eta(k)-\eta\rangle\nnum\\
&+ 2\alpha(k)\langle\hat{\DD}(k)-\tilde{\DD}(k),\eta(k)-\eta\rangle\nnum\\
&\leq \alpha(k)^2\|\DD(k)\|^2 + \|\eta(k) - \eta\|^2\nnum\\&- \|\eta(k+1) - \eta\|^2 + 2\alpha(k)\|\tilde{\DD}(k)-\DD(k)\| \|\eta(k)-\eta\|\nnum\\
&+ 2\alpha(k)\|\hat{\DD}(k)-\tilde{\DD}(k)\| \|\eta(k)-\eta\|.\label{e30}\end{align}

Now let us examine the term with $\|\hat{\DD}(k)-\tilde{\DD}(k)\|$ on the right-hand side of~\eqref{e30}. Since $x(k)\in X$ and $X$ is convex and closed, it then follows from the non-expansiveness of the projection operator $\mathbb{P}_X$, that \begin{align}\|\eta(k+1)-\eta(k)\|
&=\|\mathbb{P}_K[\eta(k)-\alpha(k)\DD(k)]-\mathbb{P}_K[\eta(k)]\|\nnum\\
&\leq\alpha(k)\|\DD(k)\|.\label{e14}\end{align} Consequently,
it follows from the Lipschitiz continuity of the game map $\nabla\Theta$ and~\eqref{e30},~\eqref{e14} that \begin{align} &\|\hat{\DD}(k)-\tilde{\DD}(k)\|\nnum\\
&\leq \sum_{i\in V}(\|\hat{\DD}^{[i]}_x(k)-\tilde{\DD}^{[i]}_x(k)\| + \|\hat{\DD}^{[i]}_{\mu}(k)-\tilde{\DD}^{[i]}_{\mu}(k)\|)\nnum\\
&\leq 2N L_{\Theta}\sum_{\tau=k-\tau_{\max}}^{k-1}\|\eta(\tau+1)-\eta(\tau)\|\nnum\\
&\leq 2N L_{\Theta}\sum_{\tau=k-\tau_{\max}}^{k-1}\alpha(\tau)\|\DD(\tau)\|.\label{e45}\end{align}

Since $K$ is compact and the image of $\nabla\Theta$ is uniformly compact, ~\eqref{e45} implies that there is $\Upsilon > 0$ such that \begin{align}&2\alpha(k)\|\hat{\DD}(k)-\tilde{\DD}(k)\|\|\eta(k)-{\eta}\|\nnum\\
&\leq 2NL_{\Theta}\sum_{\tau=k-\tau_{\max}}^{k-1}2\alpha(k)\alpha(\tau)
\|\DD(\tau)\|\|\eta(k)-{\eta}\|\nnum\\
&\leq 2N L_{\Theta}\sum_{\tau=k-\tau_{\max}}^{k-1}\big(\alpha(k)^2+\alpha(\tau)^2
\|\DD(\tau)\|^2\|\eta(k)-{\eta}\|^2\big)\nnum\\
&\leq 2N L_{\Theta}\tau_{\max}\alpha(k)^2 + \Upsilon\sum_{\tau=k-\tau_{\max}}^{k-1}\alpha(\tau)^2.\label{e46}\end{align}

Now consider the term with $\|\tilde{\DD}(k)-\DD(k)\|$ on the right-hand side in~\eqref{e30}. Recall that $K$ is compact. By the Taylor expansion, we reach that \begin{align}&\DD^{[i]}_x(k) =\nabla_{x^{[i]}}f_i(x^{[i]}(k),
\Lambda_i(k))\nnum\\
&+\sum_{\ell=1}^{m_i}\mu^{[i]}_{\ell}(k)\nabla_{x^{[i]}}
G^{[i]}_{\ell}(x^{[i]}(k),\Lambda_i(k))+O(c_i(k))\textbf{1}_{n_i}.\label{e21}\end{align}

With the above relation, we have the following for $\tilde{\DD}^{[i]}(k)-\DD^{[i]}(k)$:
\begin{align}\tilde{\DD}^{[i]}_{\mu}(k)-\DD^{[i]}_{\mu}(k) = 0,\quad
\tilde{\DD}^{[i]}_x(k)-\DD^{[i]}_x(k) = O(c_i(k))\textbf{1}_{n_i}.\label{e12}\end{align}

Notice that $\|\DD(k)\|$ is uniformly bounded. Substitute~\eqref{e46} and~\eqref{e12} into~\eqref{e30}, sum over $[0,T]$, and it renders that there is some $\Upsilon', \Upsilon'' > 0$ such that the following estimate holds: \begin{align}&2\sum_{k=0}^T\alpha(k)\langle \hat{\DD}(k), \eta(k)-\eta\rangle\nnum\\
&\leq \|\eta(0) - \eta\|^2 - \|\eta(T+1) - \eta\|^2 + \Upsilon'\sum_{k=0}^T\alpha(k)^2\nnum\\
&+ \Upsilon''\sum_{k=0}^T\alpha(k)\max_{i\in V}c_i(k).\label{e26}\end{align}

Since $\{\alpha(k)^2\}$ and $\{\alpha(k)\max_{i\in V}c_i(k)\}$ are summable, then the right-hand side of~\eqref{e26} is finite when we let $T\rightarrow+\infty$. We now show the following by contradiction: \begin{align}\liminf_{k\rightarrow+\infty}\langle \hat{\DD}(k),\eta(k)-\eta\rangle\leq0.\label{e34}\end{align}
Assume that there are $k_0\geq0$ and $\epsilon>0$ such that $\langle \hat{\DD}(k),\eta(k)-\eta\rangle \geq \epsilon$ for all $k\geq k_0$. Since $\alpha(k) > 0$ and $\{\alpha(k)\}$ is not summable, then we have the following: \begin{align}&\sum_{k=0}^{+\infty}\alpha(k)\langle \hat{\DD}(k),\eta(k)-\eta\rangle\nnum\\
&\geq \sum_{k=0}^{k_0}\alpha(k)\langle \hat{\DD}(k),\eta(k)-\eta\rangle\nnum\\
&+ \sum_{k=k_0+1}^{+\infty}\alpha(k)\langle \hat{\DD}(k),\eta(k)-\eta\rangle\nnum\\
&\geq \sum_{k=0}^{k_0}\alpha(k)\langle \hat{\DD}(k),\eta(k)-\eta\rangle
+ \epsilon\sum_{k=k_0+1}^{+\infty}\alpha(k)\nnum\\
&\geq +\infty.\nnum\end{align}
We then reach a contradiction, and thus~\eqref{e34} holds. Equivalently, the following relation holds: \begin{align*}\limsup_{k\rightarrow+\infty}\langle \hat{\DD}(k),\eta-\eta(k)\rangle\geq0.\end{align*} Since $\nabla\Theta$ is closed, the above relation implies that there is a limit point $\tilde{\eta}\in K$ of the sequence $\{\eta(k)\}$ such that the following holds: $\langle \nabla\Theta(\tilde{\eta}),\eta-\tilde{\eta}\rangle\geq0, \quad \forall\eta\in K$. Since $\nabla\Theta$ is pseudo-monotone, then we have the following relation: \begin{align}\langle \nabla\Theta(\eta),\eta-\tilde{\eta}\rangle\geq0, \quad \forall\eta\in K.\label{e22}\end{align}

We now set out to show the following by contradiction: \begin{align}\langle \nabla\Theta(\tilde{\eta}),\eta-\tilde{\eta}\rangle\geq0,\quad \forall \eta\in K.\label{e32}\end{align} Assume that there is $\hat{\eta}\in K$ such that the following holds: \begin{align}\langle \nabla\Theta(\tilde{\eta}),\hat{\eta}-\tilde{\eta}\rangle < 0.\label{e33}\end{align} Now choose $\varepsilon\in(0,1)$, and define $\eta_{\varepsilon} \triangleq \hat{\eta} + \varepsilon(\tilde{\eta}-\hat{\eta})$. Since $K$ is convex, then we have $\eta_{\varepsilon}\in K$. The following holds: \begin{align}\langle \nabla\Theta(\eta_{\varepsilon}),(1-\varepsilon)(\hat{\eta}-\tilde{\eta})\rangle = \langle \nabla\Theta(\eta_{\varepsilon}),\eta_{\varepsilon}-\tilde{\eta}\rangle\geq0,\label{e28}
\end{align} where in the inequality we use~\eqref{e22}. It follows from~\eqref{e28} that the following relation holds for any $\varepsilon\in(0,1)$: \begin{align}\langle \nabla\Theta(\eta_{\varepsilon}),\hat{\eta}-\tilde{\eta}\rangle\geq0,\quad \forall \eta\in K.\label{e29}
\end{align} Since $\nabla\Theta$ is closed, letting $\varepsilon\rightarrow1$ in~\eqref{e29} gives that: $\langle \nabla\Theta(\tilde{\eta}),\hat{\eta}-\tilde{\eta}\rangle\geq0$, which contradicts~\eqref{e33}. As a result, the relation~\eqref{e32} holds. By~\cite{Palomar.Eldar:10}, the relation~\eqref{e32} implies that $\tilde{\eta}\in\mathbb{X}_{\rm UC}$. We replay $\eta$ by $\eta(k)$ in~\eqref{e32}, and establish the following:
\begin{align}\langle \nabla\Theta(\tilde{\eta}),\eta(k)-\tilde{\eta}\rangle\geq0.\label{e31}\end{align} Since $\nabla\Theta$ is pseudo-monotone, it follows from~\eqref{e31} that\begin{align}\langle \hat{\DD}(k), \eta(k)-\tilde{\eta}\rangle \geq 0.\label{e27}\end{align} Replace $\eta$ with $\tilde{\eta}$ in~\eqref{e23}, apply~\eqref{e27}, and it renders: \begin{align}\|\eta(k+1) - \tilde{\eta}\|^2 -\|\eta(k) - \tilde{\eta}\|^2
\leq\alpha(k)^2\|\DD(k)\|^2.\label{e35}\end{align} Sum up~\eqref{e35} over $[s,k]$ and we have \begin{align}\|\eta(k) - \tilde{\eta}\|^2
\leq \|\eta(s) - \tilde{\eta}\|^2 + \sum_{\tau=s}^{k-1}\alpha(\tau)^2\|\DD(\tau)\|^2.\label{e36}\end{align}
We take the limits on $k$ first and then $s$ on both sides of~\eqref{e36}. Since $\{\alpha(k)\}$ is square summable and $\{\DD(k)\}$ is uniformly bounded, we have \begin{align}\limsup_{k\rightarrow\infty}\|\eta(k) - \tilde{\eta}\|^2
\leq \liminf_{s\rightarrow\infty}\|\eta(s) - \tilde{\eta}\|^2.\nnum\end{align} It implies that $\{\|\eta(k) - \tilde{\eta}\|\}$ converges. Since $\tilde{\eta}$ is a limit point of $\{\eta(k)\}$, $\{\eta(k)\}$ converges to $\tilde{\eta}\in\mathbb{X}_{\rm UC}$. Furthermore, we have $\tilde{x}\in\mathbb{X}_{\rm CVX}$ by (P2) in Proposition~\ref{pro5}.\oprocend

\subsection{Scenario two: strongly monotone reduced game map and known $\tau_{\max}$}

In the last section, the convergence of Algorithm~\ref{ta:algo1} is ensured under the mild assumptions: $\tau_{\max}$ is unknown and the game map $\nabla\Theta$ is merely pseudo-monotone. This set of assumptions requires the usage of diminishing step-sizes. Note that diminishing step-sizes may cause a slow convergence rate. The shortcoming can be partially addressed by choosing a constant step-size when the inequality constraints are absent and (A5) in Theorem~\ref{the1} is strengthened to the following one: \begin{assumption} The map $\nabla\Theta^r$ is strongly monotone over $X$ with constant $\rho$.\label{asm5}\end{assumption}

\begin{algorithm}[htbp]\caption{The distributed gradient-based algorithm for Scenario two} \label{ta:algo2}

\textbf{Require:} Each player in $V$ chooses the initial state $x^{[i]}(0)\in X_i$.

\textbf{Ensure:} At each $k \ge 0$, each player~$i\in V$ executes the following steps:

1. Each player $i\in V$ updates its state according to the following rule:\begin{align*} x^{[i]}(k+1) = \mathbb{P}_{X_i}[x^{[i]}(k) - \alpha D^{[i]}(k)],\end{align*} where the approximate gradient $D^{[i]}(k)$ is given by: \begin{align*}[D^{[i]}(k)]_{\ell} &= \frac{1}{2c_i(k)}
\{f_i(x^{[i]}(k)+c_i(k)e^{[\ell]}_{n_i},\Lambda_i(k))\nnum\\
&-f_i(x^{[i]}(k)-c_i(k)e^{[\ell]}_{n_i},\Lambda_i(k))\},\end{align*} for $\ell = 1,\cdots,n_i$.

2. Repeat for $k = k+1$.
\end{algorithm}

Algorithm~\ref{ta:algo2} is proposed to address Scenario two. In particular, Algorithm~\ref{ta:algo2} is similar to Algorithm~\ref{ta:algo1}, but has two distinctions. Firstly, Algorithm~\ref{ta:algo2} excludes the inequality constraints. Because the game map $\nabla \Theta$ cannot be strongly monotone due to $\HH_i$ being linear in $\mu^{[i]}$. Secondly, thanks to the strong monotonicity of $\nabla \Theta$, a constant step-size replaces diminishing step-sizes. The following theorem summarizes the convergence properties of Algorithm~\ref{ta:algo2}.

\begin{theorem} Suppose the following holds:
\begin{enumerate}
\item[(B1)] the quantity $0\leq\tau_{\max}<\frac{\rho}{4NL_{\Theta^r}}$ is known;
\item[(B2)] Assumptions~\ref{asm7} and~\ref{asm1} hold;
\item[(B3)] the step-size $\alpha\in(0,\frac{\rho-4NL_{\Theta^r}\tau_{\max}}{L^2_{\Theta^r}(2+16N^2\tau_{\max})})$;
\item[(B4)] The sequence $\{\max_{i\in V}c_i(k)\}$ is summable.
\item[(B5)] Assumption~\ref{asm5} holds.
\end{enumerate}
For any initial state $x(0)\in X$, the sequence of $\{x(k)\}$ generated by Algorithm~\ref{ta:algo2} converges to $\tilde{x}\in\mathbb{X}_{\rm CVX}$.\label{the2}
\end{theorem}

\textbf{Proof:} First of all, for any $k\geq0$ and $i$, we define $\hat{\DD}^{[i]}(k) \triangleq \nabla_{x^{[i]}}f_i(x^{[i]}(k),\{x^{[j]}(k)\}_{j\in\NN^{\rm IN}_i})$, which is the gradient evaluated at the delay-free states. So, the quantity $\hat{\DD}(k)$ is free of time delays and perturbations. Similarly, for any $k\geq0$ and $i$, we define $\tilde{\DD}^{[i]}(k) \triangleq \nabla_{x^{[i]}}f_i(x^{[i]}(k),\Lambda_i(k))$ which is the gradient evaluated at the delayed states. Then the quantity $\tilde{\DD}(k)$ is free of perturbations but subject to time delays.

We then write Algorithm~\ref{ta:algo2} in the following compact form: \begin{align}x(k+1) = \mathbb{P}_X[x(k)-\alpha D(k)].\label{e39}
\end{align} It is noticed that $\tilde{x}\in{\mathbb{X}}_{\rm CVX}$ is a fixed point of the operator $\mathbb{P}_X[\cdot-\alpha \nabla\Theta^r(\cdot)]$ for any $\alpha > 0$; i.e., it holds that $\tilde{x} = \mathbb{P}_X[\tilde{x}-\alpha \nabla\Theta^r(\tilde{x})]$. By this, we have the following relations: \begin{align}&\|x(k+1)-\tilde{x}\|^2\nnum\\
&= \|\mathbb{P}_X[x(k)-\alpha D(k)] - \mathbb{P}_X[\tilde{x}-\alpha \nabla\Theta^r(\tilde{x})]\|^2\nnum\\
&\leq \|(x(k) - \tilde{x})-\alpha(D(k)-\nabla\Theta^r(\tilde{x}))\|^2\nnum\\
&= \|x(k) - \tilde{x}\|^2-2\alpha\langle x(k) - \tilde{x}, \hat{D}(k)-\nabla\Theta^r(\tilde{x})\rangle\nnum\\
&+ \alpha^2\|(D(k)-\tilde{D}(k))+(\tilde{D}(k)-\hat{D}(k))+(\hat{D}(k)-u)\|^2\nnum\\
&-2\alpha\langle x(k) - \tilde{x}, D(k)-\tilde{D}(k)\rangle\nnum\\
&-2\alpha\langle x(k) - \tilde{x}, \tilde{D}(k)-\hat{D}(k)\rangle\nnum\\
&\leq (1-2\rho\alpha+4L_{\Theta^r}^2\alpha^2)\|x(k)-\tilde{x}\|^2\nnum\\
&+ 4\alpha^2(\|D(k)-\tilde{D}(k)\|^2+\|\tilde{D}(k)-\hat{D}(k)\|^2)\nnum\\
&+2\alpha\|x(k) - \tilde{x}\|(\|D(k)-\tilde{D}(k)\|+\|\tilde{D}(k)-\hat{D}(k)\|),\label{e42}\end{align}
where we use the non-expansiveness property of the projection operator $\mathbb{P}_X$ in the first inequality, and the strong monotonicity and the Lipschitz continuity of $\nabla\Theta^r$ in the last inequality. For the term of $\|\hat{D}(k)-\tilde{D}(k)\|$ in~\eqref{e42}, one can derive the following relations: \begin{align}&\|\hat{D}(k)-\tilde{D}(k)\| \leq  NL_{\Theta^r}\sum_{\tau=k-\tau_{\max}}^{k-1}\|x(\tau+1)-x(\tau)\|\nnum\\
&\leq  NL_{\Theta^r}\sum_{\tau=k-\tau_{\max}}^{k-1}\big(\|x(\tau+1)-\tilde{x}\|
+\|x(\tau)-\tilde{x}\|\big)\nnum\\
&\leq  2NL_{\Theta^r}\sum_{\tau=k-\tau_{\max}}^{k-1}\|x(\tau)-\tilde{x}\|.\label{e47}\end{align}

Similar to~\eqref{e12}, we have \begin{align}\tilde{D}^{[i]}(k)-D^{[i]}(k) = O(c_i(k))\textbf{1}_{n_i}.\label{e13}\end{align}

Substituting~\eqref{e47} and~\eqref{e13} into~\eqref{e42} yields: \begin{align}&\|x(k+1)-\tilde{x}\|^2\nnum\\
&\leq (1-2\rho\alpha+4L_{\Theta^r}^2\alpha^2+4NL_{\Theta^r}\alpha\tau_{\max})
\|x(k)-\tilde{x}\|^2\nnum\\
&+ (32N^2L_{\Theta^r}^2\alpha^2+4NL_{\Theta^r}\alpha)
\sum_{\tau=k-\tau_{\max}}^{k-1}\|x(\tau)-\tilde{x}\|^2\nnum\\
&+4\alpha^2(\max_{i\in V}c_i(k))^2 + 2\alpha\max_{i\in V}c_i(k)\|x(k)-\tilde{x}\|,\label{e48}\end{align}
where we use the relation $2ab \leq a^2 + b^2$.

To study the convergence of $x(k)$ in~\eqref{e48}, we define the following notation: \begin{align*}&\psi(k) \triangleq x(k) - \tilde{x},\quad\|\chi(k)\| \triangleq \max_{k-\tau_{\max}\leq \tau \leq k-1}\|\psi(\tau)\|,\nnum\\
&e(k)\triangleq 4\alpha^2(\max_{i\in V}c_i(k))^2 + 2\alpha\max_{i\in V}c_i(k)\|x(k)-\tilde{x}\|.\end{align*} Then~\eqref{e48} is rewritten as follows: \begin{align}\|\psi(k+1)\|\leq a\|\psi(k)\|+b\|\chi(k)\|+e(k).\label{e43}\end{align} where $a \triangleq (1-2\rho\alpha+4L_{\Theta^r}^2\alpha^2+4NL_{\Theta^r}\alpha\tau_{\max})$, $b \triangleq \tau_{\max}(32N^2L_{\Theta^r}^2\alpha^2+4NL_{\Theta^r}\alpha)$ and the sequence $\{e(k)\}$ is diminishing. From the recursion of~\eqref{e43}, we can derive the following relation for any pair of $k>s\geq0$: \begin{align}&\|\psi(k)\| \leq a^{k-s}\|\psi(s)\| + \sum_{\tau=s}^{k-1}a^{k-\tau}b\chi(\tau) + \sum_{\tau=s}^{k-1}e(\tau)\nnum\\
&\leq a^{k-s}\|\psi(s)\| + b\sup_{s\leq\tau\leq k}\|\chi(\tau)\|\sum_{\tau=s}^{k-1}a^{k-\tau} + \sum_{\tau=s}^{k-1}e(\tau)\nnum\\
&\leq a^{k-s}\|\psi(s)\| + \frac{b}{1-a}\sup_{s\leq\tau\leq k}\|\chi(\tau)\| + \sum_{\tau=s}^{k-1}e(\tau).\label{e44}\end{align} Recall that $\{\max_{i\in V}c_i(k)\}$ and thus $\{(\max_{i\in V}c_i(k))^2\}$ are summable. Take the limits on $k$ and $s$ at both sides of~\eqref{e44}, and it gives the following relation: \begin{align}\limsup_{k\rightarrow+\infty}\|\psi(k)\| \leq \frac{b}{1-a}\limsup_{k\rightarrow+\infty}\|\chi(k)\|.\label{e49}\end{align} On the other hand, one can see that \begin{align}\limsup_{k\rightarrow+\infty}\|\chi(k)\| = \limsup_{k\rightarrow+\infty}\|\psi(k)\|.\label{e50}\end{align}

Since $\alpha\in(0,\frac{\rho-4NL_{\Theta^r}\tau_{\max}}{L^2_{\Theta^r}(2+16N^2\tau_{\max})})$, then $\frac{a}{1-b}<1$. The combination of~\eqref{e49} and~\eqref{e50} renders that \begin{align}\limsup_{k\rightarrow+\infty}\|\psi(k)\| = 0.\label{e40}\end{align} Apparently, \begin{align}\liminf_{k\rightarrow+\infty}\|\psi(k)\| \geq 0.\label{e41}\end{align} The combination of~\eqref{e40} and~\eqref{e41} establishes the convergence of $\{x(k)\}$ to $\tilde{x}$. It completes the proof.\oprocend

\section{Discussions}\label{sec:discussion}

\subsection{Comparison of two scenarios}\label{sec:comparison}

The two proposed algorithms are complementary. Algorithm~\ref{ta:algo1} can address inequality constraints, does not need to know $\tau_{\max}$ and merely requires the game map to be pseudo monotone. It comes with the price of potentially slow convergence due to the utilization of diminishing step-sizes. In contrast, Algorithm~\ref{ta:algo2} cannot deal with inequality constraints, needs to know $\tau_{\max}$ and requires the game map to be strongly monotone. It comes with the benefit of potentially fast convergence due to the utilization of a constant step-size.

\subsection{Discussion on Assumption~\ref{asm1}}\label{sec:assumption}

The proposed algorithms rely upon Assumption~\ref{asm1}. In what follows, we will provide two sufficient conditions for Assumption~\ref{asm1}. Recall that $\sigma_{i,\min}$, $\sigma_{i,\max}$,
$\sigma_{\min}$ and $\sigma_{\max}$ are defined in Section~\ref{subsection:notations}.

\subsubsection{Global Slater vectors}

\begin{assumption} There exist $\bar{x}\in X$ and $\sigma' > 0$ such that $\|-G^{[i]}(\bar{x})\|_{\infty} \geq \sigma'$ for all $i\in V$.\label{asm3}
\end{assumption}

The vector $\bar{x}$ that satisfies Assumption~\ref{asm3} is referred to as the \emph{global} Slater vector. The existence of the global Slater vector ensures the boundedness of $\mathbb{X}_{\rm UC}(M)$ as follows.

\begin{lemma} If Assumption~\ref{asm3} holds and $G^{[i]}_{\ell}$ is convex in $x$, then for any $(\tilde{x},\tilde{\mu})\in \mathbb{X}_{\rm UC}(\real^m_{\geq0})$, it holds that $\|\tilde{\mu}^{[i]}\|_{\infty} \leq \frac{\sigma_{\max}-\sigma_{\min}}{\sigma'}$ for $i\in V$.\label{lem2}
\end{lemma}

\textbf{Proof:} Pick any $\tilde{\eta}\triangleq(\tilde{x},\tilde{\mu})\in \mathbb{X}_{\rm UC}(\real^m_{\geq0})$. It holds that \begin{align}\langle \nabla\Theta(\tilde{\eta}),\tilde{\eta}-\eta\rangle\leq0.\label{e55}\end{align}
Choose $\eta = (\bar{x}^T,0^T)^T$ in~\eqref{e55}, and we have
\begin{align}&0\geq\langle \nabla\Theta^r(\tilde{x}),\tilde{x}-\bar{x}\rangle - \sum_{i\in V}\langle \tilde{\mu}^{[i]},G^{[i]}(\tilde{x})\rangle\nnum\\
&+ \sum_{i\in V}\sum_{\ell=1}^{m_i}\tilde{\mu}_{\ell}^{[i]}\langle \nabla_{x^{[i]}} G^{[i]}_{\ell}(\tilde{x}),\tilde{x}^{[i]}-\bar{x}^{[i]}\rangle\nnum\\
&\geq \langle \nabla\Theta^r(\tilde{x}),\tilde{x}-\bar{x}\rangle - \sum_{i\in V}\langle \tilde{\mu}^{[i]},G^{[i]}(\bar{x})\rangle,\label{e56}
\end{align} where the last inequality uses $\tilde{\mu}^{[i]}_{\ell}\geq0$ and Taylor theorem and the convexity of $G^{[i]}_{\ell}$ in $x$. It follows from~\eqref{e56} and Assumption~\ref{asm3} that \begin{align}\|\tilde{\mu}^{[i]}\|_{\infty}\leq\frac{-\langle \nabla\Theta^r(\tilde{x}),\tilde{x}-\bar{x}\rangle}{\sigma'}\leq \frac{\sigma_{\max}-\sigma_{\min}}{\sigma'}.\nnum\end{align}
\oprocend


\subsubsection{Private Slater vectors}

In Lemma~\ref{lem2}, each player has to access the global information of $\bar{x}$, $\sigma'$, $\sigma_{\min}$ and $\sigma_{\max}$. In what follows, we will derive a sufficient condition which only requires private information of each player~$i$ to determine an upper bound on $\mu^{[i]}$.

\begin{assumption} For each $i\in V$, there exists $\sigma_i > 0$ and $\bar{x}_i\in X_i$ such that for any $(\tilde{x},\tilde{\mu})\in\mathbb{X}_{\rm UC}$, $\|-G^{[i]}(\bar{x}^{[i]},\tilde{x}^{[-i]})\|_{\infty} \geq \sigma_i$ holds.\label{asm2}
\end{assumption}

The vector $\bar{x}^{[i]}$ that satisfies Assumption~\ref{asm3} is referred to as the \emph{private} Slater vector. One case where Assumption~\ref{asm2} holds is that $G^{[i]}$ only depends upon $x^{[i]}$ and, for each $i\in V$, there is $\bar{x}^{[i]}\in X_i$ such that $\|-G^{[i]}(\bar{x}^{[i]})\|_{\infty} \geq \sigma_i$ holds for some $\sigma_i > 0$.

\begin{lemma} If Assumption~\ref{asm2} holds, then for any $(\tilde{x},\tilde{\mu})\in \mathbb{X}_{\rm UC}(\real^m_{\geq0})$, it holds that $\|\tilde{\mu}^{[i]}\|_{\infty} \leq \frac{\sigma_{i,\max}-\sigma_{i,\min}}{\sigma_i}$ for $i\in V$.\label{lem5}
\end{lemma}

\textbf{Proof:} Pick any $(\tilde{x},\tilde{\mu})\in\mathbb{X}_{\rm UC}(\real^m_{\geq0})$. By Assumption~\ref{asm2}, there is $\bar{x}^{[i]}\in X_i$ such that \begin{align*}\|-G^{[i]}(\bar{x}^{[i]},\tilde{x}^{[-i]})\|_{\infty} \geq \sigma.\end{align*} Since $(\tilde{x},\tilde{\mu})\in\mathbb{X}_{\rm UC}$, we then have the following: \begin{align}&\HH_i(\tilde{x}^{[i]},\tilde{x}^{[-i]},0)\leq
\HH_i(\tilde{x}^{[i]},\tilde{x}^{[-i]},\tilde{\mu}^{[i]})\nnum\\
&\leq
\HH_i(\bar{x}^{[i]},\tilde{x}^{[-i]},\tilde{\mu}^{[i]})\nnum\end{align} where in the first inequality we use $0\in M_i$ and in the second inequality we use $\bar{x}^{[i]}\in X_i$. The above relations further render the following: \begin{align}f_i(\tilde{x})\leq f_i(\bar{x}^{[i]},\tilde{x}^{[-i]}) + \langle \tilde{\mu}^{[i]}, G^{[i]}(\bar{x}^{[i]},\tilde{x}^{[-i]})\rangle.\label{e24}\end{align} Recall $\|-G^{[i]}(\bar{x}^{[i]},\tilde{x}^{[-i]})\|_{\infty} \geq \sigma_i$. Then the relation~\eqref{e24} implies the relation of $\|\tilde{\mu}^{[i]}\|_{\infty} \leq \frac{\sigma_{i,\max}-\sigma_{i,\min}}{\sigma_i}$.\oprocend

\subsection{Future directions}

The current paper imposes several assumptions: (1) the compactness of $X_i$; (2) the smoothness of component functions; (3) identical $\alpha(k)$ or $\alpha$ for all the players; (4) the convexity of component functions and constraints. In particular, the compactness of $X_i$ ensures the uniform boundedness of partial gradients. The convexity guarantees that local optimum are also globally optimal in the sense of Nash equilibrium. It is of interest to relax these assumptions.

\section{Numerical simulations}\label{sec:simulations}

In this section, we will provide a set of numerical simulations to verify the performance of our proposed algorithms. The procedure ${\rm sample}(A)$ returns a uniform sample from the set $A$.

\subsection{Algorithm~\ref{ta:algo1}}\label{sec:simulations-alg1}

Consider a power network which can be modeled as an interconnected graph $\GG_p\triangleq \{\mathcal{V}_p,\EE_p\}$ where each node $i\in \mathcal{V}_p$ represents a bus and each link in $\EE_p$ represents a branch. The buses in $\mathcal{V}_p$ are indexed by $1,\cdots,N$ and bus $1$ denotes the feeder which has a fixed voltage and flexible power injection. A lossless DC model is used to characterize the relation between power generations and loads at different buses and power flows across various branches. Each bus is connected to either a power load or supply and each load is associated with an end-user. The set of load buses is denoted by $\mathcal{V}_l\subseteq \mathcal{V}_p$ and the set of supply buses is denoted by $\mathcal{V}_s\subseteq \mathcal{V}_p$.

The maximum available power supply at bus $i\in \mathcal{V}_s$ is denoted by $S_i\geq0$, and the intended power load at bus $i\in \mathcal{V}_l$ is denoted by $L_i\geq0$. If the total power supply exceeds the total intended power load, then all the intended loads can be satisfied. Otherwise, some loads need to reduce. The reduced load of end-user~$i$ is denoted by $R_i\in[0,L_i]$. The objective function of each end-user~$i$ is given by $f_i(R_i) = c_i R_i + p(\textbf{1}^T_{|\mathcal{V}_l|}(L-R))(L_i-R_i) - u_i(L_i-R_i)$. The quantity $c_iR_i$ represents the disutility induced by load reduction $R_i$ with $c_i>0$. The scalar $p(\textbf{1}^T(L-R))$ is the charged price given the total actual load $\textbf{1}^T(L-R)$. The value $u_i(L_i-R_i)$ stands for the benefit produced by load $L_i-R_i$.

Each end-user~$i$ aims to minimize its own objective function $f_i(R_i)$ by adjusting $R_i$. Such decision making is subject to the physical constraints of the power grid. The first constraint is that the total actual load cannot exceed the maximum available supply; i.e., \begin{align}
\textbf{1}_{|\mathcal{V}_l|}^T(L-R) \leq \textbf{1}_{|\mathcal{V}_s|}^TS.\label{e52}\end{align} Another set of constraints are induced by the limitations of power flows on branches. Let $H_g\in[-1,1]^{|\EE|\times|\mathcal{V}_g|}$ (resp. $H_l\in[-1,1]^{|\EE|\times|\mathcal{V}_l|}$) be the generation (resp. load) shift factor matrix. For $H_g$, the $(i,\ell)$ entry represents the power that is distributed on line $\ell$ when $1MW$ is injected into bus $i$ and withdrawn at the reference bus. Denote by $f_e^{\max}$ the maximum capacity of branch $e$. The power flow constraint for branch $e$ can be expressed as: \begin{align}-f^{\max}\leq H_g S - H_l (L-R)\leq f^{\max},\label{e53}\end{align} where $f^{\max}\triangleq{\rm col}[f_e^{\max}]_{e\in \EE_p}$, $S\triangleq {\rm col}[S_i]_{i\in\mathcal{V}_s}$, $L\triangleq {\rm col}[L_i]_{i\in\mathcal{V}_l}$ and $R\triangleq {\rm col}[R_i]_{i\in\mathcal{V}_l}$.

The above description defines a non-cooperative game among end-users in $\mathcal{V}_l$ as follows:
\begin{align}&\min_{R_i\in [0,L_i]}c_i R_i + p(\textbf{1}_{|\mathcal{V}_l|}^T(L-R))(L_i-R_i) - u_i(L_i-R_i)\nnum\\
&{\rm s.t.}\quad \textbf{1}_{|\mathcal{V}_l|}^T(L-R) \leq \textbf{1}_{|\mathcal{V}_s|}^TS\nnum\\
&\quad\quad H_g S - H_l (L-R) - f^{\max}\leq 0\nnum\\
&\quad\quad - H_g S + H_l (L-R) -f^{\max}\leq 0.\label{e54}
\end{align}

In game~\eqref{e54}, $p$ is the pricing policy of load serving entity (LSE). This is confidential information of LSE and should not be disclosed to end-users. So the objective function $f_i$ is partially unknown to end-user~$i$. In addition, the numerical values of generation and load shift factor matrices are of national security interest and are kept confidential from the public. Therefore, the power flow constraints are unknown to end-users. Assume that all the end-users can communicate with LSE. Given $R$, LSE broadcasts the price value $p(\textbf{1}_{|\mathcal{V}_l|}^T(L-R))$, the power imbalance $\textbf{1}_{|\mathcal{V}_l|}^T(L-R) - \textbf{1}_{|\mathcal{V}_s|}^TS$, and power limit violations $H_g S - H_l (L-R) - f^{\max}$ and $- H_g S + H_l (L-R) -f^{\max}$. In this way, end-users can access the values of $f_i$ and $G^{[i]}$ in Algorithm~\ref{ta:algo1} without knowing their structures.\footnote{Many networked engineering systems operate in a hierarchical structure; e.g., Internet, power grid and transportation systems, where a central system operator is placed at the top layer and end-users are placed at the bottom layer~\cite{MC-SHL-ARC-JCD:07a,Wu.Moslehi.Bose:89}. Here we assume that LSE can communicate with all the end-users. This assumption is widely used in; e.g., network flow control~\cite{LL-99}.}

We choose $u_i(L_i-R_i) = -\frac{1}{2}a_i(L_i-R_i)^2+b_i(L_i-R_i)$ with $a_i>0$. For the pricing mechanism, we use (4) in~\cite{Bulow.Peiderer:83}; i.e., $p(w) = \beta w^{\frac{1}{\eta}}$ with $\beta>0$, and $\eta\in(0,1)$. So $f_i(R) = c_iR_i + (L_i-R_i)\beta (\textbf{1}_{|\mathcal{V}_l|}^T(L-R))^{\frac{1}{\eta}} +\frac{1}{2}a_i(L_i-R_i)^2-b_i(L_i-R_i)$. One can compute the following: \begin{align}\frac{\partial^2f_i}{\partial^2R_i}&=a_i-2\beta\frac{1}{\eta}
(\textbf{1}_{|\mathcal{V}_l|}^T(L-R))^{\frac{1}{\eta}-1}\nnum\\
&-\beta (L_i-R_i)\frac{1}{\eta}(\frac{1}{\eta}-1)(\textbf{1}_{|\mathcal{V}_l|}^T(L-R))^{\frac{1}{\eta}-2},\nnum\\
\frac{\partial^2f_i}{\partial R_i \partial R_j}&=-\beta\frac{1}{\eta}
(\textbf{1}_{|\mathcal{V}_l|}^T(L-R))^{\frac{1}{\eta}-1}\nnum\\
&-\beta (L_i-R_i)\frac{1}{\eta}(\frac{1}{\eta}-1)(\textbf{1}_{|\mathcal{V}_l|}^T(L-R))^{\frac{1}{\eta}-2}.\nnum\end{align}

\begin{figure*}[bh]
  \centering
  \includegraphics[width = \linewidth]{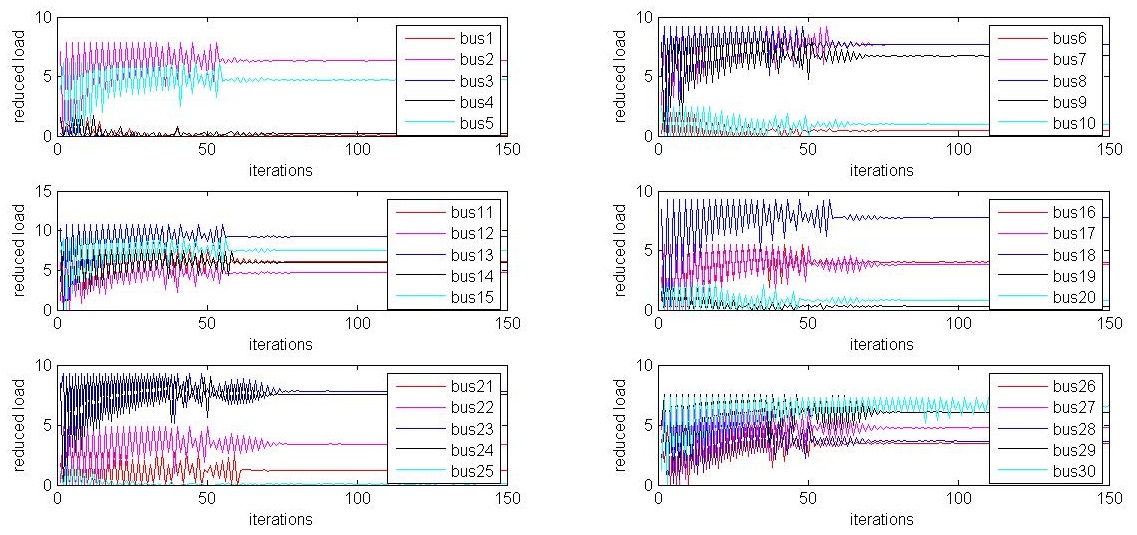}
  \caption{The estimates of Algorithm~\ref{ta:algo1} at $30$ buses for $\tau_{\max} = 10$.}\label{fig1}
\end{figure*}

Choose $a_i$ such that $a_i > (N-2)\beta\frac{1}{\eta}
(\textbf{1}_{|\mathcal{V}_l|}^TL)^{\frac{1}{\eta}-1}+(N-1)\beta L_i\frac{1}{\eta}(\frac{1}{\eta}-1)(\textbf{1}_{|\mathcal{V}_l|}^TL)^{\frac{1}{\eta}-2}>0$. Then the monotonicity property holds. Assume that there is a global Slater vector. Figure~\ref{fig1} shows the simulation results for the IEEE $30$-bus Test System~\cite{PSTCA} and the system parameters are adopted from MATPOWER~\cite{MATPOWER}. The delays at each iteration are randomly chosen from $0$ to $10$.

\subsection{Algorithm~\ref{ta:algo2}}\label{sec:simulations-alg2}

Consider $10$ players and their components functions are defined as follows: \begin{align*}f_i(x) = \frac{1}{2}((Q^{[i]})^Tx)^2 + (q^{[i]})^Tx,\end{align*} where
\begin{itemize}
\item $X_i = [\psi_i^m\;\;\psi_i^M]$ with $\psi_i^m = {\rm sample}([-10\;\;-5])$ and $\psi_i^M = {\rm sample}([5\;\;10])$;
\item $Q^{[i]}\in\real^{10}$ with $Q^{[i]}_i = 1$ and $Q^{[i]}_j = \frac{1}{10}$ for $j\neq i$;
\item $q^{[i]}\in\real^{10}$ with $q^{[i]}_j = {\rm sample}([-50\;\;50])$ for $j\in V$.
\end{itemize}

\begin{figure}[h]
  \centering
  \includegraphics[width = \linewidth]{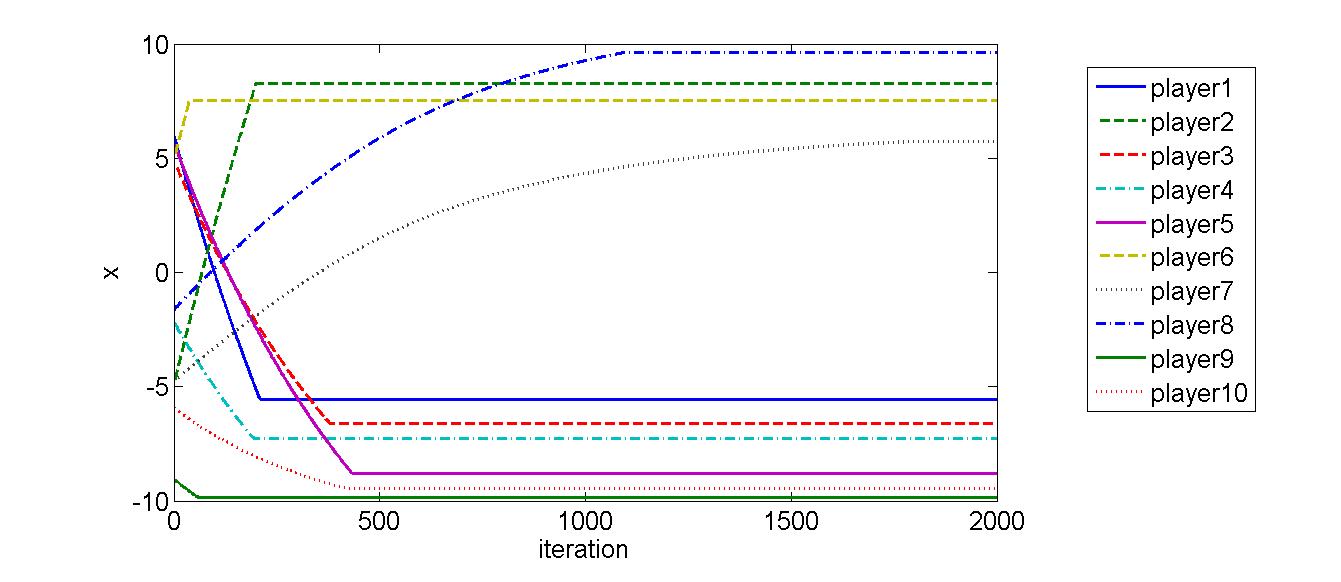}
  \caption{The estimates of Algorithm~\ref{ta:algo2} for $\tau_{\max} = 10$.}\label{fig4}
\end{figure}

One can verify that $\rho = \frac{1}{10}$ and $L_{\Theta^r} = 2$ in Theorem~\ref{the2}. By (B4), we have $\varrho_b>1$. In the simulation, we choose $\varrho_b = 3$. Figure~\ref{fig4} presents the estimate evolution of the players.

\section{Conclusions}\label{sec:conclusions}

We have studied a set of distributed robust adaptive computation algorithms for a class of generalized convex games. We have shown their asymptotic convergence to Nash equilibrium in the presence of network unreliability and the uncertainties of component functions and illustrated the algorithm performance via numerical simulations.



\bibliographystyle{plain}


\end{document}